\documentclass[10pt]{amsart}
\usepackage{amsfonts}
\usepackage{amssymb}
\usepackage[latin1]{inputenc}
\newtheorem{thm}{Theorem}[section]
\newtheorem{cor}[thm]{Corollary}
\newtheorem{lemma}[thm]{Lemma}
\newtheorem{sublemma}[thm]{Sublemma}
\newtheorem{prop}[thm]{Proposition}
\theoremstyle{definition}

\theoremstyle{remark}

\newcommand{\bbR}{\mathbb R}
\newcommand{\bbT}{\mathbb T}

\newcommand{\bbZ}{\mathbb Z}

\numberwithin{equation}{section}

\newcommand{\rank}{{\rm rank\ }}

\begin{document}

\title[Equivariant normal form for singular orbits]{Equivariant normal form for nondegenerate singular orbits of
integrable Hamiltonian systems \\
{\it Formes normales \'{e}quivariantes pour les orbites singuli\`{e}res
nond\'{e}g\'{e}n\'{e}r\'{e}es des syst\`{e}mes Hamiltoniens int\'{e}grables }}

\author{Eva Miranda}
\address{Departament d'\`{A}lgebra i Geometria, Facultat de Matem\`{a}tiques, Universitat de Barcelona,
Gran Via de les Corts Catalanes 585, 08007 Barcelona, Spain }
\email{evamiranda@ub.edu}
\thanks{ The first author is partially supported by the DGICYT project number BFM2003-03458. }

\author{Nguyen Tien Zung}
\address{Laboratoire Emile Picard, UMR 5580 CNRS, UFR MIG, Universit\'{e} Toulouse III}
\email{tienzung@picard.ups-tlse.fr}

\date{July 12, 2004}

% ----------------------------------------------------------------
\begin{abstract}
  We consider an integrable Hamiltonian system with $n$-degrees of
freedom whose first integrals are invariant under the symplectic
action of a compact Lie group $G$. We prove that the singular
Lagrangian foliation associated to this Hamiltonian system is
symplectically equivalent, in a $G$-equivariant way, to the
linearized foliation in a neighborhood of a  compact singular
non-degenerate orbit. We also show that the non-degeneracy
condition is not equivalent to the non-resonance condition for
smooth systems.

\vspace{5mm}

\noindent {\sc{ R\'{e}sum\'{e}}}. On consid\`{e}re un syst\`{e}me
hamiltonien int\'{e}grable \`{a} n degr\'{e}s de libert\'{e} et
une action symplectique d'un groupe de Lie compact $G$ qui laisse
invariantes les int\'{e}grales premi\`{e}res. On prouve que le
feuilletage lagrangien singulier attach\'{e} \`{a} ce syst\`{e}me
hamiltonien est symplectiquement \'{e}quivalent, de fa\c{c}on
$G$-\'{e}quivariante, au feuilletage linearis\'{e} dans un
voisinage d'une orbite compacte singuli\`{e}re. On d\'{e}montre
aussi que la condition de non-d\'{e}g\'{e}n\'{e}r\'{e}scence n'est
pas \'{e}quivalente \`{a} la non-r\'{e}sonance pour les
syst\`{e}mes diff\'{e}rentiables.

\end{abstract}
\maketitle

\section{Introduction}

In this paper, we are interested in the geometry of integrable
Hamiltonian systems on symplectic manifolds. When we refer to an
integrable Hamiltonian system, we mean that it is integrable in
the sense of Liouville. That is to say, the system is given by a
moment map $\bf F$ on $(M^{2n},\omega)$,
\begin{equation}
{\bf F} = (F_1,\dots,F_n): (M^{2n},\omega) \rightarrow \bbR^n
\end{equation}
whose component functions $F_i$ are functionally independent
almost everywhere and Poisson commuting ($\{F_i,F_j\} = 0$ for any
$i,j$). The Hamiltonian system considered will be called smooth or
real analytic if the corresponding moment map is so. Let $X=X_H$
be the Hamiltonian vector field associated to a given function
$H$, we say that this Hamiltonian vector field is integrable if
there exists a moment map with $F_1=H.$

In our approach to the study of the integrable Hamiltonian system,
 the original Hamiltonian function defining the system will be left aside
  and our study will be focused on
the singular Lagrangian fibration given by the level sets of the
moment map.
  This moment map
generates an infinitesimal Poisson $\bbR^n$-action on
$(M^{2n},\omega)$ via the Hamiltonian vector fields
$X_{F_1},\dots,X_{F_n}$. Denote by $O$ an orbit of this
$\bbR^n$-action. We will assume that $O$ is a closed submanifold
(i.e. compact without boundary) of $(M^{2n},\omega)$. Then it is
well-known that $O$ is diffeomorphic to a torus since the vector
fields $X_{F_1},\dots,X_{F_n}$ are complete on $O$, and $O$ is a
quotient of the Abelian group $\bbR^n$ by a cocompact subgroup.

The complete integrability of Hamiltonian systems is an old
problem.  In  the XIX century, Joseph Liouville
\cite{Liouville-1855} proved that if a Hamiltonian system has $n$
functionally independent integrals in involution then it is
integrable by quadratures.
 To our knowledge, Henri Mineur was the first who gave a complete
 description, up to symplectomorphism, of the Hamiltonian system in a neighborhood of a compact regular
 orbit of dimension $n$.
 In his papers \cite{Mineur-AA1935,Mineur-AA1936,Mineur-AA1937}, it
 is proven that under the previous assumptions, there is a
symplectomorphism $\phi$ from a neighborhood $({\mathcal
U}(O),\omega)$ of $O$ in $(M^{2n},\omega)$ to $(D^n \times \bbT^n,
\sum_1^n d\nu_i \wedge d\mu_i)$, where $(\nu_1,...,\nu_n)$ is a
coordinate system on a ball $D^n$, and $(\mu_1 (mod 1),...,\mu_n
(mod 1))$ is a periodic coordinate system  on the torus $\bbT^n$,
such that $\phi_{\ast}{\bf F}$ is a map which depends only on the
variables $\nu_1,...,\nu_n$. The functions $q_i = \phi^{\ast}
\mu_i$ on ${\mathcal U}(O)$ are called angle variables, and the
functions $p_i = \phi^{\ast} \nu_i$  are called action variables.

 Although the works of Henri Mineur date back to the thirties, the theorem stated above has been known in the literature as
Arnold-Liouville theorem.

Mineur \cite{Mineur-AA1935,Mineur-AA1937} also showed that the
action functions $p_i$ can be defined via the period integrals:
\begin{equation}
\label{eqn:MineurFormula} p_i (x) = \int_{\Gamma_i (x)} \beta
\end{equation}
Here $\beta$ is a primitive $1$-form of the symplectic form, i.e.
$d\beta = \omega$, and $\Gamma_i(x)$ for each point $x$ near $O$
is a closed curve which depends smoothly on $x$ and which lies on
the Liouville torus containing $x$. The homology classes of
$\Gamma_1(x),...,\Gamma_n(x)$ form a basis of the first homology
group of the Liouville torus. The above important formula for
finding action functions will be called Mineur's formula. It can
be used to find action functions and hence torus actions and
normalization not only near regular level sets of the moment map,
but also near singular level sets as well, see e.g.
\cite{Zung-Degenerate2000,Zung-Birkhoff2002}. In particular, in
\cite{Zung-Birkhoff2002} this Mineur's formula was used in the
proof of the existence of a local analytic Birkhoff normalization
for any analytic integrable Hamiltonian system near a singular
point.

 The above-mentioned action-angle coordinates entail a
 \lq\lq uniqueness"  for the symplectic structure and the
 regular Lagrangian fibration in a neighborhood of a compact
 orbit. In fact, they provide
  a \lq\lq linear model" in a neighborhood of a regular compact orbit.
  In the same spirit, the aim of the present paper is to
establish an analog of this classical Liouville-Mineur-Arnold
theorem for the case when the orbit $O$ is singular, i.e. is of
dimension $m = \dim O$ smaller than $n$, under a natural
nondegeneracy condition. We will show that the system can be
\lq\lq linearized" near $O$
 up to fibration-preserving symplectomorphisms. The fibration in
question is the singular Lagrangian fibration given by the moment
map. We also take into account the possible symmetries of the
system. Namely, we will show that in the case there exists a
symplectic action of a compact Lie group in a neighborhood of $O$
preserving the moment map, this linearization can be carried out
in an equivariant way.

\section{Preliminaries and statement of the main results}

\subsection{ Nondegenerate orbits}

In this paper,  the orbit $O$ is assumed to be nondegenerate. The
concept of non-degenerate orbit was introduced by Eliasson. Let us
recall what it means (see, e.g.,\cite{Eliasson-NF1990},
\cite{Zung-AL1996}). A point $x \in M^{2n}$ is called singular for
the system if its rank, i.e. the rank of $\bf F$ at $x$, is
smaller than $n$. If $\rank x = m$ then $m$ is also the dimension
of the orbit through $x$ of the local Poisson $\bbR^n$ action. If
$\rank x = 0$, i.e. $x$ is a fixed point, then the quadratic parts
$F^{(2)}_1,...,F^{(2)}_n$ of the components $F_1,...,F_n$ of the
moment map at $x$ are Poisson-commuting  and they form an Abelian
subalgebra, $\mathcal A$, of the Lie algebra $Q(2n,\bbR)$ of
homogeneous quadratic functions of $2n$ variables under the
standard Poisson bracket. Observe that the algebra $Q(2n,\bbR)$ is
isomorphic to the symplectic algebra $sp(2n,\bbR)$. A point $x$
will be called a non-degenerate fixed point if $\mathcal A$ is a
Cartan subalgebra of $Q(2n,\mathbb R)$.

More generally, when $\rank x = m \geq 0$, we may assume without
loss of generality that $dF_1 \wedge ... \wedge dF_m (x) \neq 0$,
and a local symplectic reduction near $x$ with respect to the
local free $\bbR^m$-action generated by the Hamiltonian vector
fields $X_{F_1},...,X_{F_m}$ will give us an $m$-dimensional
family of local integrable Hamiltonian systems with $n-m$ degrees
of freedom. Under this reduction, $x$ will be mapped to a fixed
point in the reduced system, and if this fixed point is
nondegenerate according to the above definition, then $x$ is
called a nondegenerate singular point of rank $m$ and corank
$(n-m)$. The orbit $O$ will be called nondegenerate if it contains
a nondegenerate singular point. In fact,  if a point in $O$ is
nondegenerate then every point of $O$ is nondegenerate because
nondegeneracy is a property which is invariant under the local
Poisson $\bbR^n$-action.

\subsection{ The Williamson type of an orbit}

According to \cite{Zung-AL1996} we will define the Williamson type
of a nondegenerate singular point $x \in O$ as a triple of
nonnegative integers $(k_e,k_h,k_f)$, where $k_e$ (resp., $k_h,
k_f$) is the number of elliptic (resp., hyperbolic, focus-focus)
components of the system at $x$.  Let us recall what $k_e$, $k_h$
and $k_f$ stand for.  When $\rank x=0$, a generic linear
combination of the linear parts of the Hamiltonian vector fields
$X_{F_1},...,X_{F_n}$ at $x$ has $k_e$ pairs of purely imaginary
eigenvalues, $k_h$ pairs of real eigenvalues, and $k_f$ quadruples
of non-real non-purely-imaginary complex eigenvalues (note that
the set of eigenvalues is symmetric with respect to the real axis
and the imaginary axis). If $\rank x\neq 0$, we can perform a
symplectic reduction first and the values of $k_e,k_h$ and $k_f$
coincide with the values of $k_e,k_h$ and $k_f$ at the point
corresponding to $x$ in the reduced space. In particular, we have
$k_e + k_h + 2k_f = n -m$. The triple $(k_e,k_h,k_f)$ is also
called the Williamson type of $O$, because it does not depend on
the choice of $x$ in $O$. When $k_h = k_f =0$, we say that the
singular orbit is of elliptic type.

\subsection{ The linear model}

We are going to introduce the linear model associated to the orbit
$O$ for a given symplectic action preserving the system. Later, we
will see that the invariants associated to the linear model are
the Williamson type of the orbit and a twisting group $\Gamma$
attached to it.

Denote by $(p_1,...,p_m)$ a linear coordinate system of a small
ball $D^m$ of dimension $m$, $(q_1 (mod 1),...,q_m (mod 1))$ a
standard periodic coordinate system of the torus $\bbT^m$, and
$(x_1,y_1,...,x_{n-m},y_{n-m})$ a linear coordinate system of a
small ball $D^{2(n-m)}$ of dimension $2(n-m)$. Consider the
manifold
\begin{equation}
\label{eqn:V} V = D^m \times \bbT^m \times D^{2(n-m)}
\end{equation}
 with the standard symplectic form $\sum dp_i
\wedge dq_i + \sum dx_j \wedge dy_j$, and the following moment
map:
\begin{equation}
\label{eqn:linearmm} ({\bf p},{\bf h}) =
(p_1,...,p_m,h_1,...,h_{n-m}): V \rightarrow \bbR^n
\end{equation}
where
\begin{equation}
\label{eqn:h_i}
\begin{array}{l}
h_i = x_i^2 + y_i^2 \ \ {\rm for} \ \ 1 \leq i \leq k_e \ ,  \\
 h_i = x_iy_i \ \ {\rm for} \ \ k_e+1 \leq i \leq k_e+k_h \ , \\
%\left.
%\begin{array}{l}
h_i = x_i y_{i+1}- x_{i+1} y_i \ \ {\rm and} \\
                            h_{i+1} = x_i y_i + x_{i+1} y_{i+1}
%\end{array}
%\right\}
\ \ {\rm for} \ \ i = k_e+k_h+ 2j-1, \ 1 \leq j \leq k_f
\end{array}
\end{equation}

Let $\Gamma$ be a group with a symplectic action $\rho(\Gamma)$ on
$V$, which preserves the moment map $({\bf p},{\bf h})$. We will
say that the action of $\Gamma$ on $V$ is {\it linear} if it
satisfies the following property:

{\it $\Gamma$ acts on the product $V = D^m \times \bbT^m \times
D^{2(n-m)}$ componentwise; the action of $\Gamma$ on $D^m$ is
trivial, its action on $\bbT^m$ is by translations (with respect
to the coordinate system $(q_1,...,q_m)$), and its action on
$D^{2(n-m)}$ is linear with respect to the coordinate system
$(x_1,y_1,...,x_{n-m},y_{n-m})$}.

Suppose now that $\Gamma$ is a finite group with a free symplectic
action $\rho(\Gamma)$ on $V$, which preserves the moment map and
which is linear. Then we can form the quotient symplectic manifold
$V/\Gamma$, with an integrable system on it given by the same
moment map as above:
\begin{equation}
\label{eqn:twistedmodel} ({\bf p},{\bf h}) =
(p_1,...,p_m,h_1,...,h_{n-m}): V/\Gamma \rightarrow \bbR^n
\end{equation}
The set $\{p_i=x_i=y_i = 0\} \subset V/\Gamma$ is a compact orbit
of Williamson type $(k_e,k_f,k_h)$ of the above system. We will
call the above system on $V/\Gamma$, together with its associated
singular Lagrangian fibration, the linear system (or linear model)
of Williamson type $(k_e,k_f,k_h)$ and twisting group $\Gamma$ (or
more precisely, twisting action $\rho(\Gamma)$). We will also say
that it is a direct model if $\Gamma$ is trivial, and a twisted
model if $\Gamma$ is nontrivial.

A symplectic action of a compact group $G$ on $V/\Gamma$ which
preserves the moment map $(p_1,...,p_m,h_1,...,h_{n-m})$ will be
called linear if it comes from a linear symplectic action of $G$
on $V$ which commutes with the action of $\Gamma$. In our case,
let $\mathcal G'$ denote the group of linear symplectic maps which
preserve the moment map then this group is abelian and therefore
this last condition is automatically satisfied. In fact ${\mathcal
G'}$ is isomorphic to $\mathbb T^m\times G_1\times G_2\times G_3$
being $G_1$ the direct product of $k_e$ special orthogonal groups
$SO(2,\Bbb R)$, $G_2$ the direct product of $k_h$ components of
type $SO(1,1,\Bbb R)$ and $G_3$ the direct product of $k_f$
components of type $\mathbb R\times SO(2,\mathbb R)$,
respectively.

Now we can formulate our main result, which is the equivariant
symplectic linearization theorem for compact nondegenerate
singular orbits of integrable Hamiltonian systems:

\begin{thm}
\label{thm:equivnormal}
 Under the above notations and assumptions, there exists a finite group $\Gamma$,
a linear system on the symplectic manifold $V/\Gamma$ given by
(\ref{eqn:V},\ref{eqn:linearmm},\ref{eqn:h_i},\ref{eqn:twistedmodel}),
and a smooth Lagrangian-fibration-preserving symplectomorphism
$\phi$ from a neighborhood of $O$ into $V/\Gamma$, which sends $O$
to the torus $\{p_i=x_i=y_i = 0\}$. The smooth symplectomorphism
$\phi$ can be chosen so that via $\phi$, the system-preserving
action of the compact group $G$ near $O$ becomes a linear
system-preserving action of $G$ on $V/\Gamma$. If the moment map
$\bf F$ is real analytic and the action of $G$ near $O$ is
analytic, then the symplectomorphism $\phi$ can also be chosen to
be real analytic. If the system depends smoothly (resp.,
analytically) on a local parameter (i.e. we have a local family of
systems), then $\phi$ can also be chosen to depend smoothly
(resp., analytically) on that parameter.
\end{thm}

Remarks.

1) In the case when $O$ is a point and $G$ is trivial, the above
theorem is due to Vey \cite{Vey-Separable1978} in the analytic
case, and Eliasson \cite{Eliasson-Thesis1984,Eliasson-NF1990} in
the smooth case. The smooth one-degree-of freedom case is due to
Colin de Verdi\`{e}re and Vey \cite{CoVe-Isochore1979}.

2) In the case when $O$ is of elliptic type and $G$ is trivial,
the above theorem is due to Dufour and Molino \cite{DuMo-AA1991}
and Eliasson \cite{Eliasson-NF1990}.

3) The analytic case with $G$ trivial of the above theorem is due
to Ito \cite{Ito-AA1991}.

4) The case with $n=2, m=1$, $O$ of hyperbolic type and $G$
trivial is due to Colin de Verdi\`{e}re and Vu Ngoc San
\cite{CoVu-2D2002}, and Curr\'{a}s-Bosch and the first author
\cite{CuMi-Linearization2002}. The symplectic uniqueness for the
linear twisted hyperbolic case is not completely established in
\cite{CuMi-Linearization2002}. A complete proof for the twisted
hyperbolic case when $n=2$ and $m=1$ was given by Curr\'{a}s-Bosch in
\cite{Cu-unpublished2001}.
 The general smooth non-elliptic
case seems to be new.

5) A topological classification for nondegenerate singular fibers
of the moment map (which  contain nondegenerate singular orbits
but are much more complicated in general) was obtained  in
\cite{Zung-AL1996}, together with the existence of partial
action-angle coordinate systems. However, the problem of
classification up to symplectomorphisms of singular orbits or
singular fibers was not considered in that paper.

6) As it was already pointed out by Colin de Verdi\`{e}re and Vu Ngoc
San, the above theorem has direct applications in the problem of
semiclassical quantization of integrable Hamiltonian systems. Of
course, it is also useful for the global study of integrable
Hamiltonian systems and their underlying symplectic manifolds.

7) As it has been shown  by the first author in \cite{miranda},
this theorem  has applications in the analogous contact
linearization problem  for completely integrable systems on
contact manifolds.

The rest of the paper is organized as follows: in section  3 we
study the case of a fixed point and give the corresponding
$G$-equivariant result.  As a by-product, we prove that the path
component of the identity of the group of symplectomorphisms
preserving the system is abelian. In section 4 we prove the
general case.
 In  the Appendix we show that the
non-degeneracy condition is not equivalent to the non-resonance
condition for smooth systems.

\section{The case of a fixed point}

In this section we consider the case when $O$ is a point and  we
prove that the symplectic action of $G$ can be linearized
symplectically in a fibration-preserving way.

This linearization result can be seen as a generalization of
Bochner's linearization theorem (\cite{bo}) in the case the action
of the group preserves additional structures: the symplectic form
and the fibration. An equivariant Darboux theorem for symplectic
actions of compact Lie groups in a neighborhood of a fixed point
was proved by Weinstein in \cite{wei}. In the case the actions
considered  are the initial action and the linear action this
equivariant Darboux theorem entails a symplectic linearization
result in a neighbourhood of a fixed point ( see for example
\cite{cha1} and \cite{wei}).

  We will linearize the action of $G$ using the
averaging methods of Bochner's linearization. These averaging
tricks will be applied to  fibration-preserving symplectomorphisms
which will be presented as the time-1-map of a Hamiltonian vector
field.

In order to linearize the action of the compact Lie group in a
fibration preserving way we will work with a linear fibration and
with the standard symplectic form. The results of Eliasson
\cite{Eliasson-Thesis1984,Eliasson-NF1990} (for smooth systems)
and Vey \cite{Vey-Separable1978} (for real analytic systems) show
that there is a fibration-preserving symplectomorphism from a
neighborhood of $O$ in $(M^{2n},\omega, {\bf F})$ to a
neighborhood of the origin of the linear system
$(\bbR^{2n},\sum_{i=1}^n dx_i \wedge dy_i, {\bf h})$, where ${\bf
h} = (h_1,...,h_n)$ is the quadratic moment map given by Formula
(\ref{eqn:h_i}). If the compact symmetry group $G$ in Theorem
\ref{thm:equivnormal} is trivial, then we are done. Suppose now
that $G$ is nontrivial. We can (and will) assume that the singular
Lagrangian fibration near $O$ is already linear. We will refer to
this Lagrangian fibration as $\mathcal F$. It remains to linearize
the action of $G$ in such a way that the fibration remains the
same. It would be interesting to adapt the proof of Eliasson for
actions of compact Lie group but unfortunately some of the steps
in his proof do not seem to  admit an equivariant version.

Let us fix some notation that we will use throughout the paper.
The vector field $X_{\Psi}$ will stand for a Hamiltonian vector
field with associated Hamiltonian function $\Psi$. We denote by
$\phi_{X}^s$ the time-$s$-map of the vector field $X$. Let $\psi$
be a local diffeomorphism $\psi: (\bbR^{2n},0) \to (\bbR^{2n},0)$.
In the sequel, we will denote by $\psi^{(1)}$ the linear part of
$\psi$ at $0$. That is to say, $\psi^{(1)}(x)=d_{0}\psi (x)$.

The group  $\mathcal G$  stands for the group of local
automorphisms preserving the system,
 ${\mathcal G}=
 \{\phi: (\bbR^{2n},0) \to (\bbR^{2n},0),{\mbox{ such that}}\quad \phi^*(\omega)=\omega,\quad
 {\bf h}\circ\phi={\bf h}\}$, and
 $\mathcal G_0$ stands for the path-component of the identity of $\mathcal
G$. We denote by $\mathfrak{g}$ the Lie algebra of germs of
Hamiltonian vector fields tangent to the fibration $\mathcal F$.

The subgroup of linear transformations contained in $\mathcal G$
is denoted by ${\mathcal G}'$. As we have observed in the
introduction ${\mathcal G}'$ is abelian.

 The goal of this section
is to prove a local linearization result (Proposition
\ref{local2}) for a given smooth action of a compact Lie group
$G$. In order to prove this result we will have to show that given
any local automorphism $\psi\in\mathcal G$ then $\psi^{(1)}\circ
\psi^{-1}$ can be presented as the time-1-flow of a Hamiltonian
vector field as it is shown in  Corollary \ref{local}.

As we will see this is, in fact, a consequence of theorem
\ref{thm:exponential} which shows that the exponential mapping
$\exp\colon \mathfrak{g}\longrightarrow \mathcal G_0$ determined
by the time-1-flow of a vector field $X\in\mathfrak{g}$ is a
surjective group homomorphism.

Before stating this theorem we need the following sublemma.
\begin{sublemma}\label{sublemma1}
The Lie algebra $\mathfrak{g}$ is abelian and for any pair of
vector fields $X_{G_1}$ and $X_{G_2}$ contained in $\mathfrak{g}$
the following formula holds
$$\phi_{X_{G_1}+X_{G_2}}^s= \phi_{X_{G_1}}^s\circ
\phi_{X_{G_2}}^s.$$
\end{sublemma}

{\it Proof}

Let $X_{G_1}$ and $X_{G_2}$ be two vector fields in
$\mathfrak{g}$.
 Since $\{G_1,G_2\}=\omega(X_{G_1},X_{G_2})$ and $X_{G_1}$ and
 $X_{G_2}$ are tangent to the Lagrangian fibration $\mathcal F$ then $\{G_1,G_2\}_{L}=0$
for any regular fiber $L$ of $\mathcal F$.  On the other hand,
since the set of regular fibers is dense  and $X_{G_1}$ and
$X_{G_2}$ are also tangent along the singular fibers,  the bracket
$\{G_1,G_2\}$ vanishes everywhere.

This implies in  turn that $[X_{G_1},X_{G_2}]=0$ and the Lie
algebra $\mathfrak{g}$ is abelian.  Therefore the flows associated
to $X_{G_1}$ and $X_{G_2}$ commute. As a consequence
$\alpha_s=\phi_{X_{G_1}}^s\circ \phi_{X_{G_2}}^s$ is a
one-parameter subgroup.

A simple computation shows that its  infinitesimal generator is
$X_{G_1}+X_{G_2}$ and this ends the proof of the sublemma.

 \hfill $\square$

Observe that given a vector field $X_G$  in $\mathfrak{g}$, its
time-s-flow $\phi_{X_G}^s$ preserves the moment map $\bf{h}$
because $X_G$ is tangent to $\mathcal F$. It also preserves the
symplectic form since it is the flow of a Hamiltonian vector
field. Finally since the vector field $X_G$ vanishes at the origin
the mapping $\phi_{X_G}^s$ fixes the origin. Therefore,
$\phi_{X_G}^s$ is contained in $\mathcal G$. In fact it is
contained in $\mathcal G_0$ since $\phi_{X_G}^0=Id$.

 We denote by $\exp\colon \mathfrak{g}\longrightarrow \mathcal G_0$
the exponential mapping defined by $exp(X_G):=\phi_{X_G}^1$ for
any $X_G\in \mathfrak{g}$. We can now state and prove the first
theorem of this section.

\begin{thm}\label{thm:exponential} The exponential $\exp\colon\mathfrak{g}\longrightarrow
\mathcal G_0$ is a surjective group homomorphism, and moreover
there is an explicit right inverse given by
$$\phi\in\mathcal G_0\longmapsto \int_0^1 X_t dt \in \mathfrak{g}$$
\noindent where $X_t\in \mathfrak{g}$ is defined by

$$X_t(R_t)=\frac{d R_t}{dt}$$

for any $C^1$ path $R_t$ contained in $\mathcal G_0$ connecting
the identity to $\phi$. \noindent

\end{thm}

{\it Proof}

The formula proved in  sublemma \ref{sublemma1} with $s=1$ shows
that the exponential is a group homomorphism. It remains to show
that it is surjective.

 Let $R_s$ be a $C^1$ path  in $\mathcal G_0$ such that
 $R_0=Id$ and
$R_1=\phi$ for a given $\phi$ $\in \mathcal G_0$.

We define the time-dependent vector field $X_t$ by the following
formula:

$$X_t(R_t)=\frac{d R_t}{dt}.$$

Observe that $X_t$ is tangent to the fibration $\mathcal F$ for
any $t$ contained in $[0,1]$ because the diffeomorphism $R_t$
preserves the fibration $\mathcal F$, $\forall t$. On the other
hand since $R_t$ is a symplectomorphism for any $t$, the vector
field $X_t$ is locally Hamiltonian. Since the symplectic manifold
considered is a neighbourhood $U$ of the origin the vector field
is indeed Hamiltonian.  Thus, $X_t$ is contained in the Lie
algebra $\mathfrak{g}$.

 Now consider $Y_t=\int_0^t X_r dr$.
 This vector field is also contained in $\mathfrak{g}$.

   We will show that $\exp Y_t= R_t$ for any  $t \in [0,1]$.
 Particularizing $t=1$, this shows that $\exp Y_1=R_1=\phi$ and
 therefore the mapping of the statement is an explicit right
 inverse of the exponential mapping and the exponential is
 surjective.

In order to show the equality $\exp Y_t= R_t,
 t \in [0,1]$. We will show that it satisfies the same
 nonautonomous differential equation

 $$X_t(\exp Y_t)=\frac{d(\exp Y_t)}{dt}$$

 \noindent and this will imply  $\exp Y_t= R_t,
 t \in [0,1]$ since the initial conditions are the same.

In fact, we will prove that,
\begin{equation}
sX_t(\phi_{Y_t}^s)=\frac{d\phi_{Y_t}^s}{dt} \label{eqn:trick}
\end{equation}
\noindent which leads to the desired result when $s=1$ since by
definition $\phi_{Y_t}^1=\exp Y_t$.

Observe that the formula we want to prove is equivalent to the
fact that the vector field  $s X_t$ is tangent to the curve
$\phi_{Y_{t+u}}^s\circ (\phi_{Y_t}^s)^{-1}$ at any point $p$.
Therefore, we can write formula \ref{eqn:trick} as,

$$sX_t=\frac{d}{du}_{\vert u=0}(\phi_{Y_{t+u}}^s(\phi_{Y_t}^s)^{-1}).$$

After differentiation with respect to $s$ the formula we want to
prove becomes:
\begin{equation}
X_t={\frac{d}{ds}\frac{d}{du}}_{\vert
u=0}(\phi_{Y_{t+u}}^s(\phi_{{Y_t}}^s)^{-1}) \label{eqn:tangent}
\end{equation}

We will first compute,
$$\frac{d}{ds}(\phi_{Y_{t+u}}^s(\phi_{Y_t}^s)^{-1})$$

 According to sublemma \ref{sublemma1} the
Lie algebra $\mathfrak{g}$ is abelian and we may write

$$(\phi_{Y_{t+u}}^s(\phi_{Y_t}^s)^{-1})=\phi_{Y_{t+u}-Y_t}^s$$

Observe that $Y_{t+u}-Y_t=\int_t^{t+u} X_r dr$.

On the one hand using the  definition of flow,

$$\frac{d}{ds}(\phi_{Y_{t+u}}^s(\phi_{Y_t}^s)^{-1})=\frac{d}{ds}(\phi_{\int_t^{t+u} X_r dr})=
\int_t^{t+u} X_r dr.$$

On the other, since

\begin{equation}
\lim_{u\rightarrow 0}\frac{\int_t^{t+u} X_r dr-u X_t}{u}=0,
\end{equation}

we can write $\int_t^{t+u} X_r dr=uX_t+o(u)$, uniformly in $t$.

Therefore,
\begin{equation}
\frac{d}{ds}(\phi_{Y_{t+u}}^s(\phi_{Y_t}^s)^{-1})=u X_t+ o(u).
\end{equation}

finally differentiating in $u$ and particularizing $u=0$ we obtain

$$\frac{d}{du}\frac{d}{ds}_{\vert u=0}(\phi_{Y_{t+u}}^s(\phi_{Y_t}^s)^{-1})_{\vert u=0}=X_t.$$

This proves formula \ref{eqn:tangent} and this ends the proof of
the theorem.

\hfill{$\square$}

 \vspace{5mm}

{\bf Remark:} Observe that this exponential mapping is not always
injective. Since a vector field $X\in\mathfrak{g}$ is a
Hamiltonian vector field tangent to $\mathcal F$,  its hamiltonian
function  is a first integral of the system given by $\bf{h}$.
Therefore $X=X_{\phi(h_1,\dots,h_n)}$ when restricted to each
connected component of the regular set of $\bf h$. Bearing this in
mind, it is easy to see that the exponential is injective if there
are only hyperbolic components (Williamson type $(0,n,0)$). If
there are elliptic or focus-focus components any vector field of
type $X=X_{2\pi k h_i}, k\in\mathbb Z$ (with $h_i$ standing for an
elliptic function or for a function $h_i$ in a focus-focus pair
$h_i,h_{i+1}$) is contained in the kernel of the exponential. In
fact, the kernel is generated by these vector fields. In
particular, this guarantees that $\exp$ is always locally
injective.

\vspace{5mm}
 The theorem above has direct applications to the linearization
 problem posed at the beginning of this section but it also tells
 us that $\mathcal{G}_0$ is abelian.

\begin{cor}\label{myabelian}

 The group ${\mathcal G_0}$ is abelian.

\end{cor}

\begin{proof}

According to Theorem \ref{thm:exponential}, the exponential
mapping is a surjective morphism of groups and according to
sublemma \ref{sublemma1} the Lie algebra $\mathfrak{g}$ is
abelian. This implies that $\mathcal G_0$ is abelian.
\end{proof}

 {\bf Remarks:}

1) One could also check that $\mathcal G_0$ is abelian using the
following: Observe that it is enough to check that any two
diffeomorphisms  $\phi_1$ and $\phi_2$ in $\mathcal G_0$ commute
on an open dense set.  We consider  the dense set $\Omega$
determined by the regular points of the fibration.
    Now consider
the  submanifold
$L_{(\delta_1,\dots,\delta_n)}=\{(x_1,y_1,\dots,x_n,y_n),\quad
x_i=\delta_iy_i, \quad \forall i\}$ with $\delta_i\in\{-1,+1\}$.
It is a Lagrangian submanifold. In the case there are no
hyperbolic components, taking  $\delta_i=1$ and $\delta_{i+1}=-1$
for the focus-focus pairs $h_i$ and $h_{i+1}$, the submanifold
$L_{(\delta_1,\dots,\delta_n)}$ is transversal to the regular
Lagrangian fibration induced by $\bf F$ on $\Omega$. So we may
apply a result of Weinstein \cite{weinstein} which ensures that
the foliation is symplectomorphic in a neighborhood of $L$ to the
foliation by fibers in $T^*(L)$ endowed with the Liouville
symplectic structure. In this way we may assume that the
symplectic form is $\omega=\sum_i dp_i\wedge dq_i$ being $q_i$
coordinates on $L$. The fibration is then determined by
${\bf{F}}=(q_1,\dots,q_n)$. Now any diffeomorphism lying in
$\mathcal G$ can be written in the form
$\phi(q_1,\dots,q_n,p_1,\dots,p_n)=(q_1,\dots,q_n,p_1+\alpha_1(q),\dots,
p_n+\alpha_n(q))$ for certain smooth  functions $\alpha_i$.
Clearly any two diffeomorphism of this form commute and therefore
this proves that $\mathcal G_0$ is abelian. In the case there are
hyperbolic components since $\phi_1$ and $\phi_2$ leave each
orthant invariant. In each orthant we may consider an appropriate
choice of $\delta_i$ for hyperbolic functions $h_i$ such that
$L_{(\delta_1,\dots,\delta_n)}$ is a transversal Lagrangian
submanifold to the fibration restricted to this orthant. And we
may repeat the argument above for $F$ restricted to each orthant
to conclude that $\mathcal G_0$ is abelian.

 As observed by
Weinstein in \cite{weinstein}, the study of local
symplectomorphism preserving the foliation by fibers in $T^*(L)$
has relevance in the study of lagrangian-foliated symplectic
manifolds.

2) Although $\mathcal G$ is also abelian for analytical systems,
it is not always abelian if we consider smooth systems as the
following example shows:

Consider $n=1$ and $h=xy$. Let $\psi$ be the smooth function:

\[\psi(x,y)=\begin{cases} {e^{-{\frac{1}{(xy)^2}}}} & \quad
x\geq 0\\ {2 e^{-{\frac{1}{(xy)^2}}}} & \quad x\leq 0\end{cases}\]

\noindent and let $\phi$ be the time-1-map of $X_{\psi}$. Then
$\phi$ does not commute with the involution $I(x,y)=(-x,-y)$.

\vspace{5mm}

Another interesting consequence of Theorem \ref{thm:exponential}
is the following result about the local automorphisms of the
linear integrable system $(\bbR^{2n},\sum_{i=1}^n dx_i \wedge
dy_i, {\bf h})$.

\vspace{5mm}

\begin{cor}\label{local}
Suppose that $\psi: (\bbR^{2n},0) \to (\bbR^{2n},0)$ is a local
symplectic diffeomorphism of $\bbR^{2n}$ which preserves the
quadratic moment map ${\bf h} = (h_1,...,h_n)$. Then,

\begin{enumerate}
\item The linear part $\psi^{(1)}$ is also a system-preserving
symplectomorphism.

\item There is a vector field contained in $\mathfrak{g}$ such
that its time-1-map is $\psi^{(1)} \circ \psi^{-1}$. Moreover, for
each vector field $X$ fulfilling this condition there is a unique
local smooth function $\Psi: (\bbR^{2n},0) \rightarrow \bbR$
vanishing at $0$ which is a first integral for the linear system
given by $\bf h$ and such that $X=X_{\Psi}$. If $\psi$ is real
analytic then $\Psi$ is also real analytic.
\end{enumerate}
\end{cor}

\begin{proof}

 We are going to construct a path connecting
$\psi$ to $\psi^{(1)}$ contained in  ${\mathcal G}=
 \{\phi: (\bbR^{2n},0) \to (\bbR^{2n},0),{\mbox{ such that}}\quad \phi^*(\omega)=\omega,\quad
 {\bf h}\circ\phi={\bf h}\}$. Given a map $\psi\in {\mathcal G}$, we
 consider

\[S_t^{\psi}(x)=\begin{cases}{\displaystyle{\frac{\psi\circ g_t}{t}}(x)} &\quad
t\in(0,2]\\ \hfill \psi^{(1)}(x) & \quad t=0
\end{cases}\] being $g_t$ the homothecy $g_t(x_1,\dots,x_n)=t(x_1,\dots,x_n)$.

 Observe that in case $\psi$ is smooth, this mapping $S_t^{\psi}$ is
smooth and depends smoothly on $t$. In case $\psi$ is real
analytic, the corresponding $S_t^{\psi}$ is also real analytic and
depends analytically on $t$.

First let us check that
 $h\circ S_t^{\psi}=h$ when $t\neq 0$. We do it component-wise.

 Let $x=(x_1,\dots,x_n)$, then
 $$h_j\circ(\frac{\psi\circ g_t}{t})(x)=\frac{(h_j\circ\psi\circ g_t)(x)}{t^2}= \frac{h_j\circ
 g_t(x)}{t^2}=h_j(x)$$ where in the first and the last equalities
 we have used the fact that each component $h_j$ of the moment map $h$ is a quadratic
 polynomial whereas  the condition
  $h\circ\psi=h$ yields the second  equality.

 Now we check that $({S_{t}^{\psi}}) ^*(\omega)=\omega$ when $t\neq 0$. Since
 $\omega=\sum dx_i\wedge dy_i$, then $g_t^*(\omega)=t^2 \omega$.
 But since $\psi$ preserves $\omega$ then $$({S_t^{\psi}}) ^*(\omega)=(\frac{\psi\circ
 g_t}{t})^*\omega=\omega$$ when $t\neq 0$.

So far we have checked the conditions $h\circ S_t^{\psi}=h$ and
$({S_{t}^{\psi}}) ^*(\omega)=\omega$ when $t\neq 0$ but since
$S_t^{\psi}$ depends smoothly on $t$ we also have that   $h\circ
{S_{0}}^{\psi}=h$ and $({S_{0}^{\psi}}) ^*(\omega)=\omega$. So, in
particular, we obtain that ${S_0^{\psi}}=\psi^{(1)}$ preserves the
moment map and the symplectic structure and therefore $\psi^{(1)}$
is also contained in ${\mathcal{G}}$. This proves the first
statement of the corollary.

In order to prove the second statement we only need to show that
$\psi^{(1)}\circ\psi^{-1}$ is contained in $\mathcal G_0$. Then,
we can apply Theorem \ref{thm:exponential} to conclude.

Consider
$$R_t= \psi^{(1)}\circ S_t^{({\psi}^{-1})}$$ with $t\in[0,1]$,
this path connects the identity to $\psi^{(1)}\circ \psi^{-1}$ and
is contained in ${\mathcal G_0}$.

Then the formula of Theorem \ref{thm:exponential} applied to this
path provides a vector field $X$ whose time-1-map is
$\psi^{(1)}\circ \psi^{-1}$ and  there exists a unique Hamiltonian
function $\Psi$ vanishing at $0$ such that $X_\Psi=X$. Since the
vector field $X_{\Psi}$ is tangent to the foliation then $\{\Psi,
h_i\}=0, \forall i$, in other words, $\Psi$ is a first integral of
the system.

\end{proof}

\vspace{5mm} In the case the action of $G$ depends on parameters
we have a parametric version of Corollary \ref{local}.

\begin{cor}\label{myparameters}

Let $D_p$ stand for a disk centered at $0$ in the parameters
$p_1,\dots,p_m$. We denote by ${\bf p}=(p_1,\dots,p_m)$. Assume
that $\psi_{\bf{p}}:  (\bbR^{2n},0) \to (\bbR^{2n},0)$ is a local
symplectic diffeomorphism of $\bbR^{2n}$ which preserves the
quadratic moment map ${\bf h}$ and which depends smoothly on the
parameters $\bf p$. Then there is a  local smooth function
$\Psi{\bf{p}}: (\bbR^{2n},0) \rightarrow \bbR$ vanishing at $0$
 depending smoothly on $\bf p$ which is a first integral for the linear system given by $\bf h$
and such that $\psi_0^{(1)} \circ \psi_{\bf{p}}^{-1}$ is the
time-1 map of the Hamiltonian vector field $X_{\Psi_{\bf{p}}}$ of
$\Psi_{\bf{p}}$. If $\psi_{\bf{p}}$ is real analytic and depends
analytically on the parameters then $\Psi_{\bf{p}}$ is also real
analytic and depends analytically on the parameters.

\end{cor}

\begin{proof}

We will apply again Theorem \ref{thm:exponential}. We consider the
path
$$M_t=\psi_0^{(1)}\circ (\psi_{g_t(\bf{p})})^{-1},$$
\noindent where $g_t({p})=(tp_1,\dots,tp_m)$.
 This path is smooth (resp. analytic) if $\psi$ is smooth (resp.
analytic) and depends analytically on $t$ and  is contained in
$\mathcal G_0$. Because of theorem \ref{thm:exponential} we can
associate
   a Hamiltonian vector field $X_{\Psi_{\bf{p}}}$ to this path such that its time-1-map
coincides with $\psi_0^{(1)}\circ (\psi_{\bf{p}})^{-1}$. By
construction, this function $\Psi_{\bf{p}}$ is smooth (or real
analytic) if $\psi_{\bf{p}}$ is smooth (or real analytic) and
depends smoothly or analytically on the parameters if
$\psi_{\bf{p}}$  does so.

\end{proof}

\vspace{5mm}

After this digression, we will prove our local linearization
result.
 By abuse of language, we will denote the local (a priori
nonlinear) action of our compact group $G$ on
$(\bbR^{2n},\sum_{i=1}^n dx_i \wedge dy_i, {\bf h})$ by $\rho$.
For each element $g \in G$, denote by $X_{\Psi(g)}$ the
Hamiltonian vector field constructed via the formula explicited in
Theorem \ref{thm:exponential} applied to the path $R_t$ explicited
in Corollary \ref{local}. The  time-1 map of this vector field is
$\rho(g)^{(1)} \circ \rho(g)^{-1}$ where $\rho(g)^{(1)}$ denotes
the linear part of $\rho(g)$. It is clear that this defines a
smooth function $\Psi(g)$.

Consider the averaging of the family of vector fields
$X_{\Psi(g)}$ over $G$ with respect to the normalized Haar measure
$d\mu$ on $G$. That is to say,

\begin{equation}
X_G(x) = \int_{G} X_{\Psi(g)}(x) d\mu,\quad x\in\Bbb R^{2n}
\label{eqn:averaging}
\end{equation}

This vector field is Hamiltonian with Hamiltonian function
$\int_{G} \Psi(g)d\mu$. Denote by $\Phi_G$ the time-1 map of this
vector field $X_G$.

\begin{prop}\label{local2}
$\Phi_G$ is a local symplectic variable transformation of
$\bbR^{2n}$ which preserves the system $(\bbR^{2n},\sum_{i=1}^n
dx_i \wedge dy_i, {\bf h})$ and under which the action of $G$
becomes linear.
\end{prop}

\begin{proof}

Since $\Phi_G$ is the time-$1$ map of vector field contained in
$\mathfrak{g}$,  $\Phi_G$ is a symplectomorphism preserving the
fibration.  Therefore it defines a local symplectic variable
transformation. Let us check that this transformation linearizes
the action of $G$.

 From the
definition of $\Phi_G$ and formula \ref{eqn:averaging},

$$\Phi_G(x)=\phi_{X_G}^1(x)=\int_G\phi_{X_{\Psi(g)}}^1(x) d\mu$$
But since, $\phi_{X_{\Psi(g)}}^1=\rho(g)^{(1)}\circ\rho(g)^{-1}$
we have,
$$ \Phi_G(x)=\int_G \rho(g)^{(1)}\circ\rho(g)^{-1}(x) d\mu$$
We proceed as in the proof of Bochner's linearization Theorem
\cite{bo},
 $$(\rho(h)^{(1)}\circ\Phi_G\circ{\rho(h)}^{-1})(x)=
\rho(h)^{(1)}\circ\int_G ( \rho(g)^{(1)} \circ \rho(g)^{-1})
({\rho(h)}^{-1}(x))d\mu.$$

Using the linearity of $\rho(h)^{(1)}$ and the fact that $\rho$
stands for an action, the expression above can be written as,
$$\int_G( \rho(h)\circ\rho(g))^{(1)} \circ
(\rho(h)\circ\rho(g))^{-1}(x) d\mu$$
 Finally this
expression equals $\Phi_G$ due to the left invariance
 property of averaging
and we have proven $\Phi_G\circ\rho(h)=\rho(h)^{(1)}\circ\Phi_G$
as we wanted.

 \end{proof}

{\bf Remark}. A fortiori, one can show that any analytic action of
a compact Lie group $G$ on $(\bbR^{2n},\sum_{i=1}^n dx_i \wedge
dy_i, {\bf h})$ must be linear (so no need to linearize). Only in
the smooth non-elliptic case the action of $G$ may be nonlinear.
And even in the smooth case, if $G$ is connected then its action
is also automatically linear.

In the case the action of $G$ depends on parameters, this
proposition and  Corollary \ref{local}  lead to its parametric
version.

\begin{prop}\label{myparameters2}

In the case the action  $\rho_{\bf p}$ depends smoothly (resp.
analytically) on parameters there exists a local symplectic
variable transformation of $\mathbb R^{2n}$, ${\Phi_{\bf{p}}}$
which preserves the system and which satisfies,

$${\Phi_{\bf{p}}}\circ\rho_{\bf p}(h)=\rho_{0}(h)^{(1)}\circ{\Phi_{\bf{p}}}$$

\end{prop}

The proposition above will be a key point in the proof of the
linearization in a neighborhood of the orbit.

\vspace{5mm}

{\bf Remark:}
 In this section we have addressed  a linearization
problem with a foliation determined by a moment map
${\bf{h}}=(h_1,\dots,h_n)$ corresponding to a non-degenerate
singularity. Thanks to the smooth linearization result of Eliasson
this moment map has very specific  component functions $h_i$ of
elliptic, hyperbolic and focus-focus type.

However, some of the results in this section do no use this
particularity and remain valid in a more general context. For
instance sublemma \ref{sublemma1} and Theorem
\ref{thm:exponential} hold for a completely integrable system
which defines a generically Lagrangian foliation and having the
origin as singular point.

Corollaries \ref{local} and \ref{myparameters} also remain valid
if we also assume that the component functions $h_i$ are
homogeneous because the path $S_t$ also preserves $\bf{h}$ under
this condition. In particular, the final linearization results
Proposition \ref{local2} and Proposition \ref{myparameters2} also
hold for foliations whose moment map has homogeneous component
functions.

\section{The general case}

Suppose now that $\dim O = m > 0$. For the moment, we will forget
about the group $G$, and try to linearize the system in a
non-equivariant way first.

First let us recall the following theorem proved by the second
author in \cite{Zung-AL1996}:

\begin{thm}
Let $(U(N),\mathcal F)$ be a nondegenerate singularity of
Williamson type $(k_e,k_h,k_f)$ of an integrable Hamiltonian
system with $n$ degrees of freedom. Then there is a natural
Hamiltonian action of a torus $\mathbb T^{n-k_h-k_f}$ which
preserves the moment map of the integrable Hamiltonian system.
This action is unique, up to automorphism of $\mathbb
T^{n-k_h-k_f}$ and it is free almost everywhere in $U(N)$.
\end{thm}

In this theorem $U(N)$ stands for a neighborhood of a leaf $N$. If
we consider an orbit instead of a leaf of rank $m$ and Williamson
type $(k_e,k_h,k_f)$ we obtain a locally free system-preserving
torus $\bbT^m$-action in a neighborhood of the orbit  $O$. In
fact, this action can be found by using either Mineur's formula
(\ref{eqn:MineurFormula}) or alternatively the flat affine
structure on  the (local) regular orbits near $O$ of the Poisson
$\bbR^{n}$-action, and the existence of $m$ non-vanishing cycles
on these orbits. Let us denote by $(p_1,...,p_m)$ an $m$-tuple of
action functions near $O$ which generates such a locally-free
$\bbT^m$-action, and denote by $X_1,...,X_m$ the $m$ corresponding
periodic (of period 1) Hamiltonian vector fields.

Of course, $O$ is an orbit of the above $\bbT^m$-action. Denote by
$\Gamma$ the isotropy group of the action of $\bbT^m$ on $O$. So
$\Gamma$ is a finite Abelian group. There is a normal finite
covering $\widetilde{{\mathcal U}(O)}$ of a tubular neighborhood
${\mathcal U}(O)$ of $O$ such that the $\bbT^m$-action on
${\mathcal U}(O)$ can be pulled back to a free $\bbT^m$-action on
$\widetilde{{\mathcal U}(O)}$. The symplectic form $\omega$, the
moment map $\bf F$ and its corresponding singular Lagrangian
fibration, and the action functions $p_1,...,p_m$ can be pulled
back to $\widetilde{{\mathcal U}(O)}$. We will use $\widetilde{}$
to denote the pull-back: for example, the pull-back of $O$ is
denoted by $\widetilde{O}$, and the pull-back of $p_1$ is denoted
by $\widetilde{p_1}$. The free action of $\Gamma$ on
$\widetilde{{\mathcal U}(O)}$ commutes with the free
$\bbT^m$-action. By cancelling out the translations
(symplectomorphisms given by the $\bbT^m$-action are called
translations), we get another action of $\Gamma$ on
$\widetilde{{\mathcal U}(O)}$ which fixes $O$. We will denote this
later action by $\rho'$.

Take a point $\widetilde{x} \in \widetilde{O}$, a local disk
$\widetilde{P}$ of dimension $(2n-m)$ which intersects
$\widetilde{O}$ transversally at $\widetilde{x}$ and which is
preserved by $\rho'$. Denote by $\hat{q}_1,...\hat{q}_m$ the
uniquely defined functions modulo 1 on $\widetilde{{\mathcal
U}(O)}$ which vanish on $\widetilde{P}$ and such that
$\widetilde{X_i} (q_i) =1$, $\widetilde{X_i} (q_j) = 0$ if $i \neq
j$. Then each local disk  $\{\widetilde{p_1} =
const.,...,\widetilde{p_m} = const.\} \cap \widetilde{P}$ near
$\widetilde{x}$ has an induced symplectic structure, induced
singular Lagrangian fibration of an integrable Hamiltonian system
with a nondegenerate fixed point,  which is invariant under the
action $\rho'$ of $\Gamma$. Applying the result of the previous
Section, i.e. Theorem \ref{thm:equivnormal} in the case with a
fixed point, compact symmetry group $\Gamma$, and parameters
$p_1,...,p_m$, we can define local functions
$\widetilde{x_1},\widetilde{y_1},...,\widetilde{x_{n-m}},\widetilde{y_{n-m}}$
on $\widetilde{P}$, such that they form a local symplectic
coordinate system on each local disk $\{\widetilde{p_1} =
const.,...,\widetilde{p_m} = const.\} \cap \widetilde{P}$, with
respect to which the induced Lagrangian fibration is linear and
the action  $\rho'$ of $\Gamma$ is linear. We extend
$\widetilde{x_1},\widetilde{y_1},...,\widetilde{x_{n-m}},\widetilde{y_{n-m}}$
to functions on $\widetilde{{\mathcal U}(O)}$ by making them
invariant under the action of $\bbT^{m}$. Define the following
symplectic form on $\widetilde{{\mathcal U}(O)}$:
\begin{equation}
\widetilde{\omega_1} = \sum \widetilde{dp_i} \wedge d \widetilde{q_i} + \sum
d\widetilde{x_i} \wedge d\widetilde{y_i}
\end{equation}
Consider the difference between $\widetilde{\omega}$ and $\widetilde{\omega_1}$:

\begin{lemma}\label{lemma4.1}  There exist functions $\widetilde{g_i}$ in a neighborhood of
$\widetilde{O}$ in $\widetilde{{\mathcal U}(O)}$, which are invariant under the
$\bbT^m$-action, and such that
\begin{equation}
\widetilde{\omega_1} - \widetilde{\omega} = \sum d \widetilde{p_i} \wedge d \widetilde{g_i} \ .
\end{equation}
\end{lemma}

{\it Proof}. By definition of $\widetilde{\omega_1}$, the vector field
$\widetilde{X_i}$ is also the Hamiltonian vector field of $\widetilde{p_i}$ with
respect to $\widetilde{\omega_1}$. Thus we have $\widetilde{X_i} \lrcorner
(\widetilde{\omega_1} - \widetilde{\omega}) =0$ and $\widetilde{X_i} \lrcorner
d(\widetilde{\omega_1} - \widetilde{\omega}) =0$. In other words,
$\widetilde{\omega_1} - \widetilde{\omega} $ is a basic $2$-form with respect to the
fibration given by the orbits of the $\bbT^m$-action, i.e. it can be viewed as a
2-form on the $(2n-m)$-dimensional space of variables
$(\widetilde{p_1},...,\widetilde{p_m},\widetilde{x_1},\widetilde{y_1},...,\widetilde{x_{n-m}},\widetilde{y_{n-m}})$.
Moreover, by construction, the restriction of $\widetilde{\omega_1}$ on each
subspace $\{\widetilde{p_1} = const.,...,\widetilde{p_m} = const.\}$ coincides with
the restriction of $\widetilde{\omega}$ on that subspace. Thus we can write
$$
\widetilde{\omega_1} - \widetilde{\omega} = \sum d \widetilde{p_i} \wedge
\widetilde{\alpha_i}
$$
where each $\widetilde{\alpha_i}$ is an 1-form in variables
$(\widetilde{p_1},...,\widetilde{p_m},\widetilde{x_1},\widetilde{y_1},...,\widetilde{x_{n-m}},\widetilde{y_{n-m}})$.
We have that $\sum d \widetilde{p_i} \wedge d\widetilde{\alpha_i} =
d\widetilde{\omega_1} - d\widetilde{\omega} =0$, which implies that the restriction
of $\widetilde{\alpha_i}$ on each subspace $\{\widetilde{p_1} =
const.,...,\widetilde{p_m} = const.\}$ is closed (hence exact), so we can write
$\widetilde{\alpha_i}$ as
$$
\widetilde{\alpha_i} = d \widetilde{\beta_i} + \sum \widetilde{c_{ij}} d
\widetilde{p_j}
$$
where $\widetilde{\beta_i}$ and $\widetilde{c_{ij}}$ are functions of variables
$(\widetilde{p_1},...,\widetilde{p_m},\widetilde{x_1},\widetilde{y_1},...,\widetilde{x_{n-m}},\widetilde{y_{n-m}})$.
Thus
$$
\widetilde{\omega_1} - \widetilde{\omega} = \sum_{i <j} (\widetilde{c_{ij}}-
\widetilde{c_{ji}}) d\widetilde{p_i} \wedge d\widetilde{p_j} + \sum d\widetilde{p_i}
\wedge d\widetilde{\beta_i}
$$
Since $d(\widetilde{\omega_1} - \widetilde{\omega}) = 0$, the 2-form $\sum_{i <j}
(\widetilde{c_{ij}}- \widetilde{c_{ji}}) d\widetilde{p_i} \wedge d\widetilde{p_j}$
is closed (hence exact), and the functions $\widetilde{c_{ij}}- \widetilde{c_{ji}}$
are independent of the variables
$(\widetilde{x_1},\widetilde{y_1},...,\widetilde{x_{n-m}},\widetilde{y_{n-m}})$.
Thus we can write
$$
\sum_{i <j} (\widetilde{c_{ij}}- \widetilde{c_{ji}})
d\widetilde{p_i} \wedge d\widetilde{p_j} = d (\sum d\widetilde{p_i} \wedge
\widetilde{\gamma_i})
$$
where $\widetilde{\gamma_i}$ are functions of variables
$(\widetilde{p_1},...,\widetilde{p_m})$. Now put $\widetilde{g_i} =
\widetilde{\beta_i} + \widetilde{\gamma_i}$. \hfill $\square$

Consider the $\widetilde{g_i}$ given  by the lemma and define
\begin{equation}
\widetilde{q_i} =\hat{q_i} - \widetilde{g_i} .
\end{equation}
 Then with respect to the coordinate
system $(\widetilde{p_i},\widetilde{q_i},\widetilde{x_j},\widetilde{y_j})$, the
symplectic form $\widetilde{\omega}$ has the standard form, the singular Lagrangian
fibration is linear, and the free action of $\Gamma$ is also linear.

\vspace{3mm}
 {\bf Remark.} There is another proof of theorem \ref{thm:equivnormal} in the general case
  with trivial $G$ which
 does not use lemma \ref{lemma4.1}. It goes as follows. Assume that we have constructed the system of coordinates
$$\widetilde p_1,\hat q_1, \widetilde x_1,\widetilde y_1,\dots
,\widetilde p_m,\hat q_m \widetilde x_{n-m},\widetilde y_{n-m}$$
as before.
 Let $\mathcal D_R$ be the symplectic
distribution $$\mathcal D_R=<\frac{\partial}{\partial \widetilde
p_1},X_1,\dots,\frac{\partial}{\partial \widetilde p_m},X_m>$$ and
let $\mathcal D_S$ be the distribution symplectically orthogonal
to $\mathcal D_R$. Thus, we can write $$\widetilde w=\widetilde
\omega_R+\widetilde \omega_S.$$ It is easy to check that $D_S$ is
an involutive distribution. Denote by $\mathcal N_{S}^p$ the
integral manifold of $D_S$ through the point $p$ then
$\widetilde\omega_{S}\vert p=\widetilde\omega_{\mathcal
N_S^p}\vert p$. We can apply Theorem \ref{thm:equivnormal} to each
$\mathcal N_{S}^p$ to obtain a new system of coordinates
$\widehat{x_1},\widehat{y_1},\dots,\widehat{x_{n-m}},\widehat{y_{n-m}}$
in a neighborhood of the origin such that $$\widetilde
\omega_{S}=\sum d\widehat x_i\wedge d\widehat y_i.$$
 Finally, in
the system of coordinates $\widetilde p_1,\hat q_1, \widehat
x_1,\widehat y_1,\dots ,\widetilde p_m,\hat q_m \widehat
x_{n-m},\widehat y_{n-m}$, the symplectic form $\widetilde\omega$
has the standard form, the singular Lagrangian fibration is
linear, and the free action of $\Gamma$ is also linear.

 \vspace{5mm}

Thus we have shown that the original singular Lagrangian fibration
near $O$ is symplectically equivalent to a linear model (direct if
$\Gamma$ is trivial and twisted if $\Gamma$ is nontrivial).

As was mentioned in the introduction, the group of all linear
automorphisms (i.e. linear moment map preserving
symplectomorphisms) of the linear direct model of Williamson type
$(k_e,k_h,k_f)$ is isomorphic to the following Abelian group:
$$\bbT^m \times \bbT^{k_e} \times (\bbR \times \bbZ/2\bbZ)^{k_h} \times (\bbR \times \bbT^1)^{k_f} .$$
In particular, $\Gamma$ is necessarily a subgroup of $(\bbZ/2\bbZ)^{k_h}$. (It comes
from involutions of hyperbolic components, and it does not mix the components).

\vspace{5mm}

{\bf Remark:} In \cite{marle1} and \cite{marle2}  Marle
establishes a model for a Hamiltonian action of a compact Lie
group in a neighborhood of an orbit. This result was obtained
independently by Guillemin and Sternberg in \cite{gs1}. In our
construction a Hamiltonian action of an $m$-dimensional torus
preserving the fibration determined by the moment map comes into
the scene. In fact the linearization result we have just proved
can be understood as a generalization of Guillemin-Sternberg-Marle
theorem in the case the group considered is $\mathbb T^m$. It
gives a linear model for a Hamiltonian action of a torus
preserving an additional structure: a singular Lagrangian
fibration.

\vspace{5mm}

In order to prove Theorem \ref{thm:equivnormal} it remains to
consider the case  when the compact symmetry group $G$ is
nontrivial.

In the case there exists an action of a non-trivial group $G$ on
$V/\Gamma$ we have an induced action of $G$ on $V$ the following
theorem shows that this action can be linearized.

\begin{thm}\label{linorbit}

Let $G$ be a compact Lie group preserving the system $(D^m \times
\bbT^m \times D^{2(n-m)},\sum_{i=1}^m dp_i\wedge
dq_i+\sum_{i=1}^{n-m} dx_i \wedge dy_i, {\bf F})$ then there
exists  $\Phi_G$ a diffeomorphism defined in a tubular
neighborhood of the orbit $L=\mathbb T^m$ which preserves the
system $(D^m \times \bbT^m \times D^{2(n-m)},\sum_{i=1}^m
dp_i\wedge dq_i+\sum_{i=1}^{n-m} dx_i \wedge dy_i, {\bf F})$ and
under which the action of $G$ becomes linear.

\end{thm}

\begin{proof}

After shrinking the original neighborhood if necessary, we may
assume without loss of generality that we are considering  a
$G$-invariant neighborhood of $L$. First of all, let us express in
local coordinates how the action looks like. We denote by $\rho$
the action of $G$. For convenience, we use the simplifying
notation $p=(p_1,\dots,p_m)$ and $(x,y)=(x_1,y_1,\dots,
x_{n-m},y_{n-m})$.
 Since $G$
preserves the system, in particular $\rho$ preserves $p$ and sends
$\frac{\partial}{\partial q_i}$ to $\frac{\partial}{\partial
q_i}$. After all these considerations, for each $h\in G$ the
diffeomorphism $\rho(h)$ can be written as,

$$\rho(h)(p,q_1,\dots,q_m,x,y)=(p,q_1+g_1^h(p,x,y),\dots,q_m+g_m^h(p,x,y),\alpha^h(x,y,p))$$

\noindent where the functions $g_i^h$ and $\alpha^h$ are
constrained by more conditions given by the preservation of the
system. Before considering these constraints, it will be most
convenient to simplify the expression of $\alpha^h$ first. This
will be done using the local linearization theorem with parameters
(Proposition \ref{myparameters2}).

In order to do that, we restrict our attention to the induced
mapping,

$$\overline\rho(h)(p,x,y)=(p,\alpha^h(p,x,y))$$

\noindent and we consider the family of diffeomorphisms
$\overline\rho(h)_p:D^{2(n-m)}\longrightarrow D^{2(n-m)}$ defined
as follows,
$$\overline\rho(h)_p(x,y)=\alpha^h(p,x,y).$$

We may look at $p=(p_1,\dots,p_m)$ as parameters. For each $p$ the
mapping $\overline\rho(h)_p(x,y)$ induces an action of $G$ on the
disk $D^{2(n-m)}$ which preserves the induced system
$(D^{2(n-m)},\sum_{i=1}^n dx_i \wedge dy_i, {\bf h})$. Observe
that the preservation of the induced system implies, in
particular, that the action fixes the origin.

According to Proposition \ref{myparameters2} we can linearize the
action $\overline\rho(h)_p$ in such a way that it is taken to the
parametric-free linear action $\overline\rho(h)_0^{(1)}$. We can
extend trivially the diffeomorphism  $\Phi_p$ in the disk provided
by Proposition \ref{myparameters2} to a diffeomorphism $\Psi$ in
the whole neighborhood considered, simply by declaring,
$\Psi(p,q_1,\dots, q_m,x,y)=(p,q_1,\dots,q_m,\Phi_p(x,y))$. This
diffeomorphism does not preserve (in general) the symplectic
structure $\omega$, let $\omega_1=\Psi^*(\omega)$, we can apply
lemma \ref{lemma4.1} to $\omega$ and $\omega_1$. Consider the
$\widetilde{g_i}$ given  by lemma \ref{lemma4.1} and define
\begin{equation}
\widetilde{q_i} ={q_i} - \widetilde{g_i} .
\end{equation}
 Then after this change of coordinates $\omega_1$ is taken to
 $\omega$. For the sake of simplicity we will keep on using the
 notation $q_i$ for the new coordinates.
 After this
linearization in the $(x,y)$-direction the initial expression of
$\rho(h)$ looks like,

$$\rho(h)(p,q_1,\dots,q_m,x,y)=(p,q_1+g_1^h(p,x,y),\dots,q_m+g_m^h(p,x,y),
\overline\rho(h)_0^{(1)}(x,y)),$$

Since the action preserves the symplectic form $\sum_{i=1}^m
dp_i\wedge dq_i+\sum_{i=1}^{n-m} dx_i \wedge dy_i$ we conclude
that the functions $g_i^h$ do not depend on $(x,y)$ and so far
just depend on the parameters $(p_1,\dots, p_m)$.

That is,
$$\rho(h)(p,q_1,\dots,q_m,x,y)=(p,q_1+g_1^h(p),\dots,q_m+g_m^h(p),
\overline \rho(h)_0^{(1)}(x,y)),$$

 Observe that if we prove that these functions $g_i^h$ do
not depend on $p$ then we will be done because then the induced
action on $\mathbb T^m$ will be performed by translations. And, in
all, the action will be linear.

Consider $\mathcal H=\{\rho(h), h\in G\}$, we are going to prove
that this group is abelian.

We have to check that $\rho(h_1)\circ
\rho(h_2)=\rho(h_2)\circ\rho(h_1)$

We compute $$\rho(h_1)\circ \rho(h_2)(p,q_1,\dots,q_m,x,y)=$$
$$(p,q_1+g_1^{h_2}(p)+g_1^{h_1}(p)
,\dots,q_m+g_m^{h_2}(p)+g_m^{h_1}(p),
\overline\rho(h_1)_0^{(1)}\circ \overline\rho(h_2)_0^{(1)}(x,y))$$

\noindent on the other hand,

$$\rho(h_2)\circ
\rho(h_1)(p,q_1,\dots,q_m,x,y)=$$
$$(p,q_1+g_1^{h_1}(p)+g_1^{h_2}(p)
,\dots,q_m+g_m^{h_1}(p)+g_m^{h_2}(p),
\overline\rho(h_2)_0^{(1)}\circ \overline\rho(h_1)_0^{(1)}(x,y))$$

Clearly, the first $2m$ components coincide. As for the $2(n-m)$
last components, we can use the fact that $\mathcal G'$, the group
of linear transformations preserving the fibration and the
symplectic form is Abelian, as we pointed out in the introduction.

So far we know that the group $\mathcal H$ is abelian. It is also
compact, therefore it is a direct product of a torus $\mathbb T^r$
with finite groups of type $\mathbb Z/n \mathbb Z$. We are going
to check that for each $\rho(h)\in\mathcal H$ the functions
$g_i^h$ do not depend on $p$. It is enough to check it for
$\rho(h)$ in one of the components $\mathbb Z/n \mathbb Z$ and
$\mathbb T^r$. So we distinguish two cases,

\begin{itemize}

\item $\rho(h)$ belongs to $\mathbb Z/n
\mathbb Z$.

Then  $\rho(h)^n=Id$  this condition yields, $n g_i^{h}(p)=2\pi
m_i(p)$, $m_i(p)\in \mathbb Z$ for all $1\leq i\leq m$. Since
$m_i(p)$ is a continuous function taking values in $\mathbb Z$ it
is a constant function $m_i$. Thus, $g_i^{h}(p)=\frac{2\pi
m_i}{n}$ and $g_i^h$ does not depend on p.

\item $\rho(h)$ belongs to $\mathbb T^r$.
 We can consider  a sequence $\rho(h_n)$  lying
on the torus which belong to a finite group $\mathbb Z/k_n \mathbb
Z$ and which converge to $\rho(h)$. For each of these points
$\rho(h_n)$ we can apply the same reasoning as before to obtain,
$g_i^{h_n}(p)=\frac{2\pi m_i}{k_n}$.

Now for each $n$, the diffeomorphism $\rho(h_n)$ does not depend
on $p$, we may write this condition as,

$$\frac{\partial \rho(h_n)}{\partial p_i}=0,\quad 1\leq i\leq m$$

 Now since the action is smooth we can take limits in this
expression to obtain that
$$\frac{\partial\rho(h)}{\partial p_i}=0, \quad 1\leq i\leq m$$
and finally $g_i^{h}(p)$ does not depend on $p$.
\end{itemize}

And this ends the proof of the theorem.
\end{proof}

This linearizes the action of $G$ on $V$.  After considering the
quotient with the action of $\Gamma$ this theorem yields Theorem
\ref{thm:equivnormal}.

\section{Appendix: Nonresonance versus nondegeneracy}

As pointed out by Ito \cite{Ito-Birkhoff1989,Ito-AA1991} (see also
\cite{Zung-Birkhoff2002}), in the real analytic case, the
nondegeneracy condition explained in the introduction this paper
is essentially equivalent to the nonresonance condition. However,
in the smooth case, this is no longer true: smooth integrable
systems which are nonresonant at a singular point can be very
degenerate at that point at the same time. In particular, we have:

\begin{prop}
\label{prop:NoBirkhoff} Let $\gamma_1,...,\gamma_n$ be any
$n$-tuple of positive numbers which are linearly independent over
$\bbZ$, $n \geq m+ 2, m \geq 0$. Then, there is a smooth
integrable Hamiltonian function $H$ in a neighborhood of the
(elliptic invariant) torus $\bbT^m_0 = \bbT^m \times \{0\} \times
\{0\}$ in the standard symplectic space $(\bbT^m \times \bbR^m
\times \bbR^{2(n-m)}, \omega_0 = \sum_{i=1}^m dp_i \wedge dq_i +
\sum_{i=m+1}^n dx_i \wedge dy_i)$, such that $H = \frac{1}{\pi}
\sum_{i=1}^m \gamma_i p_i + \sum_{i=m+1}^n \gamma_i (x_i^2 +
y_i^2)  \ + $ higher order terms at $\bbT^m_0$, and such that $H$
does not admit a $C^1$-differentiable local Birkhoff normalization
near $\bbT^m_0$ (i.e. the corresponding Lagrangian fibration
cannot be linearized). Moreover, this integrable Hamiltonian
function $H$ can be chosen so that it does not admit a non-trivial
symplectic ${\mathbb S}^1$-symmetry near $\bbT^m_0$.
\end{prop}

In the above proposition, by integrability of $H$ we mean the
existence of a smooth moment map $(F_1,...,F_n)$ from a
neighborhood of $\bbT^m_0 = \bbT^m \times \{0\} \times \{0\}$  in
$ \bbT^m \times \bbR^m \times {\mathbb R}^{2(n-m)}$ to ${\mathbb
R}^n$, with $F_1 =H$, such that $\{F_i,F_j\}= 0$, and $dF_1 \wedge
... \wedge dF_n \neq 0$ almost everywhere. By  ``higher order
terms'' in $H$ we mean terms which are at least quadratic in
variables $p_i$, or cubic in variables $p_i,x_i,y_i$. According to
Sard's theorem (about the set of singular values) and Liouville's
theorem, almost all common level sets of such a moment map are
Liouville tori. The condition that $\gamma_1,...,\gamma_n$ are
independent over $\bbZ$ means that $\bbT^m \times \{0\} \times
\{0\} $ is an nonresonant invariant elliptic torus of the
Hamiltonian $H$. Recall that if there is a differentiable Birkhoff
normal form, then (since we are in the elliptic case), the system
also admits a Hamiltonian $\bbT^n$-symmetry near the elliptic
singular orbit. Thus, if the system does not admit a non-trivial
${\mathbb S}^1$-symmetry near $\bbT^m_0$, then of course it cannot
admit a differentiable Birkhoff normal form.

The proof of Proposition \ref{prop:NoBirkhoff} is inspired by what
happens to generic perturbations of integrable systems: resonant
tori that break up and give way to smaller-dimensional invariant
tori, homoclinic orbits, diffusion, etc. Usually this breaking up
of resonant tori leads to a chaotic behavior of the system, see,
e.g., \cite{Sevryuk-Invariant1998}. In order to prove Proposition
\ref{prop:NoBirkhoff}, we will construct an {\it integrable
perturbation} of the quadratic Hamiltonian ${1 \over \pi}
\sum_{i=1}^m \gamma_i p_i + \sum_{i=m+1}^n \gamma_i (x_i^2 +
y_i^2) \ + $ higher order terms at $\bbT^m_0$ in such a way that
there are also invariant tori arbitrarily
close to $\bbT^m_0$ that break up. \\

First let us consider the case with $m=0$ (i.e. a fixed point).
Our construction of $H$ in this case consists of two steps. \\

{\it Step 1. Creation of resonant regions.}

Choose a smooth function $Q(I_1,...,I_n)$ of $n$ variables
$I_1,...,I_n$ with the following properties:

a) $Q(0) = 0$, and the linear part of $Q$ at $0$ is $\sum \gamma_i
I_i$.

b) There is a series of disjoint small open balls $U_k$ in
$\bbR^n_+ = \{(I_1,...,I_n) \in \bbR^n, \ I_1 > 0,..., I_n > 0\}$,
which tend to $0$ (in Hausdorff topology) as $k \rightarrow
\infty$, such that  we have
$$Q(I_1,...,I_n) = \sum \gamma^k_i I_i \ \ \forall \ (I_1,...,I_n) \in U_{k} \ , $$
where  $\gamma^k_i$ are rational numbers such that $\lim_{k \to
\infty} \gamma^k_i = \gamma_i$.

Of course, such a function exists, and it can be chosen to be
arbitrarily close to $\sum \gamma^k_i I_i$ in
$C^{\infty}$-topology. Now put
$$H_1 = Q(x_1^2 + y_1^2,...,x_n^2 + y_n^2) .$$

Then for this integrable Hamiltonian function $H_1$, there are
open regions $V_k \subset {\mathbb R}^{2n}$ arbitrarily close to
$0$ in ${\mathbb R}^{2n}$ which are filled by resonant tori in
which the Hamiltonian flow of $H_1$ is periodic. These regions
$V_k$ are preimages of the open sets $U_k$ chosen above under the
moment map $(I_1,...,I_n)$.
\\

{\it Step 2. Creation of hyperbolic singularities.}

We will modify $H_1$ inside each open subset $V_k$ by a
$C^{\infty}$-small function which is flat on the boundary of
$V_k$, in such a way that after the modification our Hamiltonian
function remains integrable inside $V_k$ but admits an hyperbolic
singularity there. Since $H_1$ has periodic flow in $V_k$ for each
$k$, we can create a common model and then put it to each $V_k$
after necessary rescalings. The model can be done for a
2-dimensional system depending on $n-1$ parameters, and then take
a direct product of it with ${\mathbb T}^{n-1}$. In the
2-dimensional case, it is obvious how to change a regular function
on $D^1 \times {\mathbb S}^1$ into a function with a hyperbolic
singularity by a $C^{\infty}$-small perturbation. After the above
modifications, we obtain a new smooth Hamiltonian function $H_2$,
which is $C^{\infty}$-close to $H_1$, which coincides with $H_1$
outside  the union of $V_k$, and which is still smoothly
integrable (though smooth first integrals for $H_2$ will
necessarily be very degenerate at $0$). Note that, by
construction, the quadratic part of $H_2$ at $0$ is $\sum_{i=1}^n
\gamma_i (x_i^2 + y_i^2)$.

Since $H_2$ has hyperbolic singularities arbitrarily near $0$, it
cannot admit a differentiable Birkhoff normal form in a
neighborhood of $0$, for simple topological reasons concerning the
associated Liouville foliation. It is also easy to see that $H_2$
cannot admit an ${\mathbb S}^1$ symmetry near $0$: if there is a
symplectic ${\mathbb S}^1$-action in a neighborhood of $0$ which
preserves $H_2$, then this action must also preserve the
hyperbolic periodic orbits of $H_2$ in the resonant regions $V_k$.
This, in turn, implies that there is a natural number $N$ such
that $N \gamma^k_i \in \bbZ \quad , \forall k,i$, which is
impossible by our construction.
\\

Let us now consider the case $m \geq 1$. For simplicity, we will
assume that $m=1$. (The case $m > 1$ is absolutely similar). We
will repeat the above two steps to create hyperbolic
singularities, but with $n$ replaced by $n-1$ (recall that by
hypothesis $n \geq m+ 2$ so $n-1\geq 2$). The regions $V_k$  now
lie in $R^{2n-2}$, and $H_2 = \sum_{i=2}^n \gamma_i I_i + ...$ By
choosing our open sets $U_k$, we can assume that there is a
function $F = F(x_2^2 + y_2^2, ..., x_n^2 + y_n^2)$ on
$\bbR^{2n-2}$, which is flat at $0$, and such that
$$F = {\pi \over k}(x_n^2 + y_n^2) \ \ {\rm in} \ \ V_k$$
Denote by $\varphi$ the time-1 map of the Hamiltonian vector field
$X_F$ of $ F$ on $\bbR^{2n-2}$. Observe that  $\varphi$ is a
symplectomorphism formally equivalent to the identity map at $0$,
and that in each region $V_k$ the map $\varphi$ generates a
nontrivial $\bbZ/k\bbZ$ symmetry (i.e. the $k$-times iteration of
$\varphi$ in $V_k$ is the identity, and the lower iterations are
not). Now we can construct $H_2$ in such a way that in each region
$V_k$ it is also invariant under the $\bbZ/k\bbZ$ symmetry
generated by $\varphi$. Then, since $H_2 = H_1$ outside of the
sets $V_k$, and $\{H_1,F\} = 0$ by construction,  the map
$\varphi$ preserves $H_2$ everywhere in $\bbR^{2n-2}$.

We now construct our symplectic manifold using suspension. More
precisely, consider the free component-wise symplectic action of
$\bbZ$ on $\bbR^2 \times \bbR^{2n-2}$, where the action of $\bbZ$
on $\bbR^{2n-2}$ is generated by $\varphi$, and its action on
$\bbR^2$ (with coordinates $(p,q)$ and symplectic form $dp \wedge
dq$) is generated  by the shift $(p,q) \mapsto (p,q+ 1)$. Then the
suspension of $\varphi$ is the quotient of $\bbR^{2} \times
\bbR^{2n-2}$ by this $\bbZ$-action. Denote this quotient by $M$
and denote the projection  by $\pi$.  Let $\overline V_k$ be
 $\pi(\mathbb R^2\times V_k)$. Observe that the function $H_2$ is
$\varphi$-invariant and both functions, $H_2$ and $p$, are
invariant by the shift therefore $H_3 = {\gamma_1 \over \pi}p +
H_2$ can be projected to $M$. We denote  by ${\overline H}_3$ the
projection of $H_3$. Since $H_3$ is integrable on $\bbR^2 \times
\bbR^{2n-2}$ and the action defining the suspension is symplectic,
the function ${\overline H}_3$ defines an integrable Hamiltonian
system on $M$. We denote by $L_{k}$ an orbit of $X_{{\overline
H}_3}$ through an hyperbolic singularity of ${\overline H}_2$ in
the region $\overline V_k$ and denote by $L$ be the orbit of
$X_{{\overline H}_3}$ through the origin $O$. Observe that, by
construction, the orbits $L$ and $L_k$ are circles.
 Let us see that there is no
symplectic ${\mathbb S}^1$-action preserving ${\overline H}_3$ in
a neighborhood of $L$.
 Assume there existed one, then the orbits
of $X_{{\overline H}_3}$ would also be preserved by this action.
Given  an action of a Lie group $\phi:G\times M\longrightarrow M$,
we use the standard notation $G_x$ for the isotropy group at the
point $x$. We are going to use the Slice theorem to reach a
contradiction.
 \vspace{4mm}

 According to  the Slice Theorem \cite{Pa1960} for proper group actions
 there would exist a slice for the action through $x\in L$. Take $x=O$,
 we denote by $S$ the slice through the origin. From now on, we are going to
 consider a  neighborhood of $L$ invariant by this ${\mathbb S}^1$-action.
 Then, there exists a $k_0$ such that for all $k\geq k_0$,  the orbit $L_k$
 is fully contained in this neighborhood.
Since
 $S$ is a slice,  the
orbit $L_k$
  is  transverse
 to $S$  at each  $p\in S$ and   if $S\cap L_{k}$ is not empty then it consists of a finite number of points.
 Further, by construction, $S\cap L_{k}$ consists exactly of $k$
points $\{p_1,\dots, p_k\}$ for $k\geq k_1$. Those points lie in
an orbit for the action, therefore for each $p_i$ we can consider
an element
 $g_i\in {\mathbb S}^1$ such that $\phi(g_i,p_i)=p_1$ and all the
 $g_i$ obtained in this way
 are different. Following \cite{Pa1960}, if $S$ is a slice though $x$ then for $p\in S$ and $g\in G$ the condition
 $\phi(g,p)\in S$ implies  $g\in G_x$. Therefore all  the
 elements $g_i$ are contained in ${\mathbb S}_0^1$, the isotropy
 group at the origin. By construction, $k$ tends to infinity as we
 are approaching the origin, and  ${\mathbb S}_0^1$ is a compact group  containing an
 infinity of elements, therefore  ${\mathbb S}_0^1=\mathbb S^1$.
 This yields a contradiction because then the orbit through the
 origin would be reduced to a point. Therefore there exists no
  symplectic $\mathbb S^1$-action preserving $\overline H_3$.

This ends the proof of the proposition. \hfill{$\square$}
\vspace{5mm}

{\bf Acknowledgements} The authors wish to express their gratitude
to the referee for the many suggestions done. We feel that these
enlightening remarks have improved considerably the initial
version of this paper. In fact, it was following the referee's
suggestions that Theorem \ref{thm:exponential} was finally
included in the paper.

 The first author also wants to thank Carlos
Curr\'{a}s-Bosch for his remarks and Ignasi Mundet i Riera for drawing
her attention to the fact that the path used in the proof of
Corollary \ref{local} does also preserve any fibration with
homogeneous component functions.

% ----------------------------------------------------------------
\bibliographystyle{amsplain}
%\bibliography{zung}

%\end{document}

\providecommand{\bysame}{\leavevmode\hbox
to3em{\hrulefill}\thinspace}

\end{document}